\documentclass[12pt]{article}
\usepackage{amsmath}
\usepackage{graphicx}
\usepackage{enumerate}
\usepackage{natbib}
\usepackage{xcolor}
\usepackage{bm, mathtools, booktabs, adjustbox}
\usepackage{hyperref}
\usepackage{url} 
\allowdisplaybreaks
\newtheorem{theorem}{Theorem}
\newtheorem{Proposition}{Proposition}

\newtheorem{remark}{Remark}

\newcommand{\blind}{1}

\begin{document}

\def\spacingset#1{\renewcommand{\baselinestretch}%
{#1}\small\normalsize} \spacingset{1}


\if1\blind
{
  \title{\bf Impact of existence and nonexistence of pivot on the coverage of empirical best linear prediction intervals for small areas}
  \author{Yuting Chen\thanks{
Yuting Chen initiated this research while affiliated with the Joint Program in Survey Methodology at the University of Maryland.}\\
  Department of Mathematics \& Statistics, Eastern Kentucky University\\
    and \\
    Masayo Y. Hirose \\
    Institute of Mathematics for Industry, Kyushu University\\
    and \\
    Partha Lahiri\\
    Joint Program in Survey Methodology, University of Maryland}
  \maketitle
} \fi

\if0\blind
{
  \bigskip
  \bigskip
  \bigskip
  \begin{center}
    {\LARGE\bf Impact of existence and nonexistence of pivot on the coverage of empirical best linear prediction intervals for small areas}
\end{center}
  \medskip
} \fi

\bigskip
\begin{abstract}
We advance parametric bootstrap theory to construct highly efficient empirical best prediction intervals for small area means, achieving a coverage error rate of $O(m^{-3/2})$, where $m$ is the number of areas modeled by a linear mixed normal model.  For a general mixed effect model with random effects from a known but non-normal distribution (with unknown hyperparameters), we show analytically that the empirical best linear (EBL) prediction interval maintains the same coverage error order, assuming the existence of a pivot for standardized random effects with known hyperparameters.
Recognizing challenges in proving pivot existence, we develop a simple moment-based method to claim non-existence of pivot.  We find that without a pivot, the parametric bootstrap EBL prediction interval fails to reach the desired $O(m^{-3/2})$ coverage error.  We obtain a surprising result that the order $O(m^{-1})$ term is always positive under certain conditions indicating possible overcoverage of the existing parametric bootstrap EBL prediction interval.  In general, we analytically show for the first time that the coverage problem can be corrected by adopting a suitably devised double parametric bootstrap.
Our Monte Carlo simulations show that our proposed single bootstrap method  performs reasonably well when compared to rival methods.
\end{abstract}

\noindent%
\vfill

\newpage
\spacingset{1.9} 
\section{Introduction}\label{sec:1}
Small area estimation plays a vital role in survey applications, where public and private agencies depend on accurate statistical inference. While previous research has extensively addressed point prediction and associated mean squared prediction error (MSPE), interval estimation has typically been limited to specific cases, such as the linear mixed normal model. In this paper, we propose parametric bootstrap methods for constructing empirical best prediction intervals under a general area-level model. We analytically demonstrate that the proposed methods achieve high efficiency in terms of coverage error, apply to a much broader class of models than existing approaches, and avoid the need for complex analytical derivations.

The following two level model, commonly referred to as the area level model, has been extensively used in small area applications. 

\noindent\textbf{The area-level model.} For $i = 1, \cdots, m$,

Level 1 (Sampling model): $y_i\,| \,\theta_i \overset{\mathrm{ind}}{\sim} \mathcal{N}(\theta_i, D_i)$;

Level 2 (Linking model): $\theta_i \overset{\mathrm{ind}}{\sim} \mathcal{G}(\bm{x_i^\prime \beta}, A, \bm\phi)$.

\noindent Here, $m$ represents the number of small areas. The level 1 model is used to account for the sampling distribution of the direct estimates $y_i$, which are weighted or unweighted averages of observations from small area $i$. The level 2 model links the true small area means $\theta_i$ to a vector of $p$ known auxiliary variables $\bm{x_i} = (x_{i1}, \cdots, x_{ip})^\prime$, often obtained from various administrative records. The  coefficient vector $\bm{\beta} \in \mathcal{R}^p$ and model variance parameter $A$ of the linking model are generally unknown and are estimated from the available data. As in other papers on the area-level model, the sampling variances $D_i$ are assumed to be known. In practice, $D_i$'s are estimated using a smoothing technique such as the ones given in \cite{fay1979estimates}, \cite{otto1995sampling}, \cite{ha2014methods} and \cite{hawala2018variance}.

The classical area-level model introduced by \cite{fay1979estimates} assumes normality at both levels. The normality assumption at Level 1 is often justified by the Central Limit Theorem, as $y_i$ typically represents a sample summary \citep{rao2015small, jiang2022goodness}. However, the normality assumption at Level 2 is more challenging to justify. Therefore, in this work, we consider that
$\mathcal{G}$ is a fully parametric distribution, not necessarily normal, with mean $\text{E}(\theta_i) = \bm{x_i^\prime \beta}$, variance $\text{Var}(\theta_i) = A \geq 0$, and any additional parameters $\bm{\phi}$. 
For instance,
\cite{datta1995robust} assumed that $\mathcal{G}$ is a scale mixture of normal distribution. Later, \cite{bell2006using} and \cite{xie2007estimation} applied a $t$ distribution as a specific case of the scale mixture of normals at level 2, with mean $\bm{x_i^\prime \beta}$, variance $A$ and degrees of freedom $\phi$, to mitigate the influence of outliers. \cite{fabrizi2010robust} introduced two robust area-level models: the first assumes that $\theta_i$ follows an exponential power distribution with a mean $\bm{x_i^\prime \beta}$, variance $A$ and shape parameter $\phi$; the second model assumes a skewed exponential power distribution on level 2, with a mean $\bm{x_i^\prime \beta}$, variance $A$, shape parameter $\psi$ and the skewness parameter $\lambda$. Thus, in this case, the parameter vector $\bm{\phi}$ can be defined as $(\psi, \lambda)^\prime$. 
Similarly, \cite{jiang2022goodness} considered that $\theta_i$ follows a skewed normal distribution with mean \( \bm{x_i^\prime \beta} \), variance \( A \), and skewness parameter $\phi$.

The two-level model can be rewritten as the following simple linear mixed model
\begin{equation}
    y_i = \theta_i +e_i = \bm{x_i^\prime\beta} + u_i + e_i, \; i = 1, \cdots, m,
    \label{AreaModel}
\end{equation}
where the random effect $u_i \overset{\mathrm{iid}}{\sim} \mathcal{G}(0, A, \bm\phi)$ and the sampling error $e_i \overset{\mathrm{ind}}{\sim}\mathcal{N}(0, D_i)$ are independent. When $\mathcal{G}$ is different from normal distribution, the best predictor (BP) of $\theta_i$, $\tilde{\theta_i}^{\rm{BP}} = \text{E}(\theta_i|y_i)$ may not have a closed form. Instead, the best linear predictor (BLP) of $\theta_i$ always has the explicit form as below: 
\begin{equation}
\Tilde{\theta_i}^{\rm{BLP}} = (1-B_i)y_i + B_i \bm{x_i^\prime \beta},
\label{BLP}
\end{equation}
where $B_i = D_i/(A+D_i)$, and the BLP minimizes the mean squared prediction error (MSPE) among all linear predictors of $\theta_i$. 
Note that having a closed-form expression not only facilitates theoretical analysis but also simplifies practical implementation.
The variance of the prediction error $(\theta_i - \Tilde{\theta_i}^{\rm{BLP}})$ is denoted by
$g_{1i}(A) \coloneq \text{Var}(\theta_i - \Tilde{\theta_i}^{\rm{BLP}}) = AD_i/(A+D_i)$.  If $A$ is known, one can obtain the standard weighted least squares estimator of $\bm \beta$, denoting as $\Tilde{\bm\beta}(A)$. Replacing $\bm \beta$ in \eqref{BLP} by $\tilde{\bm\beta} = \tilde{\bm\beta}(A)$, a best linear unbiased prediction (BLUP) estimator of $\theta_i$ is given by 
\begin{equation}
\tilde{\theta_i}^{\text{BLUP}} = (1-B_i)y_i + B_i \bm{x}_i^\prime \tilde{\bm\beta},
\label{BLUP}
\end{equation}
which does not require normality assumptions typically used in small area estimation models. In practice, it is common that both $\bm \beta$ and $A$ are unknown and need to be estimated based on data. After plugging in their estimators, an empirical best linear unbiased predictor (EBLUP) of $\theta_i$ is given by:
\begin{equation}
\hat{\theta}_i^{\rm{EBLUP}} = (1-\hat{B}_i)y_i + \hat{B}_i \bm{x}_i^\prime \hat{\bm\beta},
\label{EBLUP}
\end{equation}
where $\hat{B}_i = D_i/(\hat{A} + D_i)$ and $\hat{\bm\beta} = \Tilde{\bm\beta}(\hat{A})$, and $\hat{A}$ is a consistent estimator of $A$ for large $m$. After plugging in the estimator $\hat{A}$, we have $g_{1i}(\hat{A}) = \hat{A}D_i/(\hat{A} + D_i)$. To simplify the notation, throughout the remainder of this paper, we will use $g_{1i}$ and $\hat{g}_{1i}$ to denote $g_{1i}(A)$ and $g_{1i}(\hat{A})$, respectively. Point prediction using EBLUP and the associated MSPE estimation have been studied extensively. See \cite{rao2015small} and \cite{jiang2007linear} for a detailed discussion on this subject.  

In this paper, we consider prediction interval approximates of small area means $\theta_i$. A prediction interval of $\theta_i$, denoted by $I_i$, is called a $100(1-\alpha)\%$ interval for $\theta_i$ if $\mathrm{P}(\theta_i \in I_i | \bm{\beta}, A, \bm{\phi}) = 1 - \alpha$, for any fixed $\bm{\beta} \in \mathcal{R}^p$, $A \in \mathcal{R}^{+}$ and $\bm{\phi} \in \Theta_{\bm{\phi}}$, the parameter space of $\bm{\phi}$. The probability $\mathrm{P}$ is with respect to the area-level model. There are several options for constructing interval estimates for $\theta_i$. The prediction interval based on only the Level 1 model for the observed data is given by $y_i \pm z_{\alpha/2}\sqrt{D_i}$, where $z_{\alpha/2}$ is the $100(1-\alpha/2)$th standard normal percentile. While the coverage probability for this direct interval is $1 - \alpha$, it is not efficient when $D_i$ is large as in the case of small area estimation. \cite{hall2006parametric} considered the interval based the regression synthetic estimator: $(\bm{x}_i^\prime\hat{\bm\beta}+b_{\hat{\alpha}_l^s}\sqrt{\hat{A}}, \bm{x}_i^\prime\hat{\bm\beta} + b_{\hat{\alpha}_u^s}\sqrt{\hat{A}})$, where $b_{\hat{\alpha}_l^s}$ and $b_{\hat{\alpha}_u^s}$ are obtained using a parametric bootstrap method (described in detail in a later section). However, this interval is synthetic in the sense that it is constructed using synthetic regression estimator of $\theta_i$ and its associated uncertainty measure, which are not area specific in the outcome variable and hence may cause larger length of the confidence interval.

It is of importance to combine both levels of the model in the interval estimation. A effective approach is to use empirical best methodology. We call an interval empirical best (EB) prediction interval if it is based on empirical best predictor of $\theta_i$. For a special case of the area-level model that $\theta_i$ follows a normal distribution, \cite{cox1975prediction} initiated the exact empirical Bayes interval: $\hat{\theta}_i^{\rm{EB}} \pm z_{\alpha/2}\sqrt{\hat{g}_{1i}}$, where $\hat{\theta}_i^{\rm{EB}}$ is the empirical Bayes estimator of $\theta_i$. Although the Cox interval always has smaller length than that of the direct interval, its coverage error is of the order $O(m^{-1})$, not accurate enough in most small area applications. \cite{yoshimori2014second} improved the cox-type empirical Bayes interval by using a carefully devised adjusted residual maximum likelihood (ARML) estimator of $A$. Their interval has a coverage error of order $O(m^{-3/2})$. Additionally, they analytically showed that their interval always produces shorter length than the corresponding direct interval. However, the properties of both the ARML estimator of $A$ and the associated interval have not yet been explored for cases involving non-normally distributed random effects.

A function of the data and parameters is called to be a pivot if its distribution does not depend on any unknown quantity \citep{shao2008mathematical, hall2013bootstrap}. When we have a linear mixed normal model on \eqref{AreaModel} ,  
$(\theta_i - \Tilde{\theta}_i^{\rm{BLP}})/\sqrt{g_{1i}}$ is a standard normal pivot. The traditional method of interval estimation for $\theta_i$ is of the form $\rm{EBLUP} \pm z_{\alpha/2}\sqrt{\text{mspe}}$, where mspe is an estimate of the true mean squared prediction error (MSPE) of the EBLUP. Unfortunately, $(\theta_i - \hat{\theta}_i^{\text{EBLUP}})/\sqrt{\hat{g}_{1i}}$ is not a pivot and this traditional approach produces too short or too long intervals. The coverage error of such interval is of the order $O(m^{-1})$, not accurate enough in most small area applications.  Recognizing that $(\theta_i -\hat{\theta}_i^{\text{EBLUP}})/\sqrt{\hat{g}_{1i}}$ does not follow a standard normal distribution, \cite{chatterjee2008parametric} and \cite{li2010adjusted} developed a parametric bootstrap method to approximate its distribution and obtained a EB prediction interval for $\theta_i$ in linear mixed normal models. They showed such interval has the coverage error of the order $O(m^{-3/2})$. However, the property remains unknown for non-normal distributed random effects. 

 In this paper, one main aim is to bring out the virtues of pivoting or rescaling, which can decrease the dependence of our statistics on unknown parameters and yield improved prediction interval approximation for small areas under the general model \eqref{AreaModel}. Analogous to \cite{chatterjee2008parametric}, we propose parametric bootstrap methods to approximate the distribution of a suitably centered and scaled EBLUP under the general model \eqref{AreaModel}, and apply it to construct a $100(1-\alpha)\%$  prediction interval estimation of $\theta_i$. 
 In fact, our theory does not require normality at either level of the area level model.
 
 Here, we define an interval based on the EBLUP of $\theta_i$ as an Empirical Best Linear (EBL) prediction interval.
 Specifically, we introduce these two key quantities: one is the centered and scaled EBLUP following the $F_{1i}$ distribution, which can be expressed as below:
 
\begin{equation*}
    H_i(\bm{\hat{\beta}}, \hat{A}) \coloneqq (\theta_i -\hat{\theta}_i^{\rm{EBLUP}})/\sqrt{\hat{g}_{1i}} \sim F_{1i},
\end{equation*}

and the other is based on BLP which has the $F_{2i}$ distribution and can be expressed as 
\begin{equation*}
H_i(\bm{\beta}, A) \coloneqq (\theta_i - \Tilde{\theta}_i^{\rm{BLP}})/\sqrt{g_{1i}} \sim F_{2i}.
\end{equation*}
If $F_{2i}$ does not depend on any unknown parameters, $H_i(\bm{\beta}, A)$ can be referred as a pivot; otherwise, $H_i(\bm{\beta}, A)$ is not a pivot and  we rewrite $F_{2i}$ as $F_{2i}(\bm{\nu})$, where $\bm{\nu} = (\nu_1, \cdots, \nu_s)^\prime$ is the unknown parameter vector determining the distribution of $H_i(\bm{\beta}, A)$. We define $q_{\alpha}$ as the $\alpha$-level quantile of the distribution $F_{2i}(\bm{\nu})$, if $\mathrm{P}(H_i(\bm{\beta}, A) \leq q_{\alpha}|\;\bm{\nu}) = \alpha$ for any fixed $\bm{\nu}$.

In Section \ref{sec:2}, we introduce a list of notation and regularity conditions used throughout the paper. In Section \ref{sec:3}, we propose a single parametric bootstrap EBL prediction interval for small area means in the context of an area-level model, where the random effects follow a general (non-normal) known distribution, except possibly for unknown hyperparameters. We analytically demonstrate that the coverage error order of the EBL prediction interval remains the same as in \cite{chatterjee2008parametric}, even when relaxing the normality of random effects through the existence of a pivot for suitably standardized random effects when hyperparameters are known. However, when a pivot does not exist, the EBL prediction interval fails to achieve the desired coverage error order, i.e. $O(m^{-3/2})$. Surprisingly, we find that the order $O(m^{-1})$ term is always positive under certain conditions, indicating potential overcoverage in the current single parametric bootstrap EBL prediction interval. Recognizing the challenge of showing existence of a pivot, we develop a simple moment-based method to claim non-existence of pivot. In Section \ref{sec:4}, we propose a double parametric bootstrap method under a general area level model and, for the first time, analytically show that this approach can correct the coverage problem. In Section \ref{sec:5}, we compare our proposed EBL prediction interval methods with the direct method, other traditional approaches, and the parametric bootstrap prediction interval proposed by \cite{hall2006parametric}.
In Section \ref{sec:6}, we demonstrate our proposed methods using real life data from the 1989 Small Area Income and Poverty Estimates (SAIPE) program.

\section{A list of notations and regularity conditions}\label{sec:2}
We use the following notations throughout the paper:
\begin{description}
    \item  $\bm{Y} = (y_1, \cdots, y_m)^\prime$, a $m\times 1$ column vector of direct estimates;
    \item $\bm{X}^\prime = (\bm{x}_1, \cdots, \bm{x}_m)$, a $p\times m$ known matrix of rank $p$;
    \item $\bm{\Sigma} = {\rm diag}(A + D_1,\cdots,A + D_m)$, a $m \times m$ diagonal matrix;
    \item $\tilde{\bm{\beta}} = (\bm{X}^\prime \bm{\Sigma} ^{-1}\bm{X})^{-1}\bm{X}^\prime \bm{\Sigma}^{-1}\bm{Y}$, weighted least square estimator of $\bm{\beta}$ with known $A$;
    \item $\bm{P} = \bm{\Sigma}^{-1} - \bm{\Sigma}^{-1}\bm{X}(\bm{X}^\prime \bm{\Sigma}^{-1}\bm{X})^{-1}\bm{X}^\prime \bm{\Sigma}^{-1}$;
    \item $\partial F_{2i}(x; \tilde{\bm{\nu}}) / \partial \bm{\nu} = \partial F_{2i}(x; \bm{\nu})/\partial \bm{\nu} |_{\bm{\nu} = \tilde{\bm{\nu}}}$, derivative with respect to $\bm{\nu}$ evaluated at $\Tilde{\bm{\nu}}$.
\end{description}

 We assume following regularity conditions for proving various results presented in this paper:

\begin{itemize}
\item[R1:] The rank of $\bm{X}$, $\text{rank}(\bm{X}) = p$, is bounded for a large $m$; 
\item[R2:] $0 < \inf_{i \ge 1} D_i \leq \sup_{i \ge 1} D_i < \infty$, $A \in (0, \infty)$ and the true $\bm{\phi} \in \Theta_{\bm{\phi}}^0$, the interior of $\Theta_{\bm{\phi}}$;
\item[R3:] $\text{sup}_{i\ge 1} h_{ii} \equiv \bm{x}_i^\prime(\bm{X}^\prime \bm{X})^{-1}\bm{x}_i = O(m^{-1})$;
\item[R4:] $\text{E}|u_i|^{8+\delta} < \infty$ for $0<\delta<1$;
\item[R5:] The distribution function $F_{2i}(\cdot)$ is three times continuously differentiable with respect to ($\cdot$), and third derivative is uniformly bounded.
When $H_i(\bm{\beta}, A)$ is not a pivot, $F_{2i}(\cdot;\bm{\nu})$ is three times continuously differentiable with respect to the parameter vector ${\bm \nu}$, and its third derivative 
is uniformly bounded; 
\item[R6:] 
When having a non-pivot $H_i(\bm{\beta}, A)$, we assume the estimator $\bm{\hat{\nu}}$ satisfies:
\begin{eqnarray*}
        \text{E}[(\hat{\bm{\nu}} - \bm{\nu})^\prime \partial F_{2i}(x; \bm{\nu})/\partial \hat{\bm{\nu}}] & = & O(m^{-1}), \\
        \text{E}\left[(\hat{\bm{\nu}} - \bm{\nu})^\prime \Bigl\{\partial^2 F_{2i}(x; \bm{\nu})/\partial\hat{\bm{\nu}}\partial\hat{\bm{\nu}^\prime}\Bigr\}(\hat{\bm{\nu}} - \bm{\nu})\right] & = & O(m^{-1}),\\
        \text{E}\lVert(\hat{\bm{\nu}} - \bm{\nu})\rVert^3 = \sum_{i, j, k}\text{E}[(\hat{\nu}_i-\nu_i)(\hat{\nu}_j-\nu_j)(\hat{\nu}_k-\nu_k)] & = & o(m^{-1}),
\end{eqnarray*}
where the expectation is taken at the true $\bm{\nu}$.
For a special case that $\bm{\nu} = A$, using the method of moment from \cite{prasad1990estimation} to estimate $A$, \cite{lahiri1995robust} showed that $\text{E}(\hat{A}_{\rm{PR}} - A) = o(m^{-1})$.
Moreover, they provided the Lemma C.1 that under the regularity conditions R1 - R5, $\text{E}|\hat{A}_{\rm{PR}} - A
|^{2q} = O(m^{-q})$ for any $q$ satisfying $0 < q \leq 2+\delta^\prime$ with $0 < \delta^\prime < \frac{1}{4}\delta$.
\end{itemize}

\section{Single parametric bootstrap}\label{sec:3}
For the remainder of this paper, without further explicit mention, we will use $\Tilde{\theta}_i$ to represent the BLP of $\theta_i$, and $\hat{\theta}_i$ to denote the EBLUP of $\theta_i$. The single parametric bootstrap method has been widely studied for its simplicity in obtaining prediction intervals directly from the bootstrap histogram. Ideally, a prediction interval of $\theta_i$ can be constructed based on the distribution of $(\theta_i - \hat{\theta}_i)/\sqrt{\hat{g}_{1i}}$, although this distribution function $F_{1i}$ is complex and difficult to approximate analytically. In this paper, we firstly follow the procedure introduced by \cite{chatterjee2008parametric} and \cite{li2010adjusted} to provide a bootstrap approximation of $F_{1i}$ by using a parametric bootstrap method. The implementation is straightforward, following these steps:

\begin{itemize}
    \item[1.] Conditionally on the data $(\bm{X, Y})$, draw $\theta_i^*$ for $ i=1,\cdots, m$, independently from the distribution $\mathcal{G}(\bm{x_{i}}^\prime\hat{\bm \beta}, \hat{A}, \hat{\bm\phi})$; 
    \item[2.] Given $\theta_i^*$, draw $y_i^*$ from the distribution $\mathcal{N}(\theta_i^*, D_i)$; 
    \item[3.] Construct the bootstrap estimators $\hat{\bm \beta}^*$, $\hat{A}^*$ using the data $(\bm{X, Y^*})$ and then obtain $\hat{\theta}_i^* = (1-\hat{B}^*_i)y_i + \hat{B}^*_i \bm{x_i^\prime \hat{\beta}^*}$ and $\hat{g}_{1i}^* = \hat{A}^*D_i/(\hat{A}^* + D_i)$;
    \item[4.] Calibrate on $\alpha$ using the bootstrap method. Let $(q_{\alpha_l}, q_{\alpha_u}) = (q_{\hat{\alpha}_l^s}, q_{\hat{\alpha}_u^s})$ be the solution of the following equations:
\begin{equation}\label{sb}
\begin{split}
        \mathrm{P}^*(H_i(\hat{\bm\beta}^*, \hat{A}^*) \leq q_{\alpha_u} \,| \, \hat{\bm \beta}, \hat{A}, \hat{\bm\phi}) &= 1 - \alpha/2\\
    \mathrm{P}^*(H_i(\hat{\bm\beta}^*, \hat{A}^*) \leq q_{\alpha_l} \,| \, \hat{\bm \beta}, \hat{A}, \hat{\bm\phi}) &= \alpha/2,
\end{split}
\end{equation}
where $H_i(\hat{\bm\beta}^*, \hat{A}^*) = (\theta_i^* -\hat{\theta}^*_i)/\sqrt{\hat{g}^*_{1i}}$;

    \item[5.] The single bootstrap calibrated prediction interval is constructed by: $$\hat{I}_i^s = \left(\hat{\theta}_i+q_{\hat{\alpha}_l^s}\sqrt{\hat{g}_{1i}} \,, \, \hat{\theta}_i+q_{\hat{\alpha}_u^s}\sqrt{\hat{g}_{1i}}\right).$$
\end{itemize}

One of our main results from the algorithm above is that when we relax the normality of the random effects by the existence of a pivot $H_i(\bm{\beta}, A)$, the prediction interval $\hat{I}_i^s$ obtained above still has a high degree of coverage accuracy. That is, it brings the coverage error down to $O(m^{-3/2})$.

\begin{theorem}
Under regularity conditions 
R1-R5,
for a preassigned $\alpha \in (0, 1)$ and arbitrary $i = 1, \cdots, m$, when $H_i(\bm{\beta}, A)$ is a pivot,  we have
\begin{equation}
\mathrm{P}(\theta_i \in \hat{I}_i^s) = 
\mathrm{P}\left(\hat{\theta}_i+q_{\hat{\alpha}_l^s}\sqrt{\hat{g}_{1i}} \leq \theta_i \leq \hat{\theta}_i+q_{\hat{\alpha}_u^s}\sqrt{\hat{g}_{1i}}\right) = 1 - \alpha + O(m^{-3/2}),
\end{equation}
where $q_{\hat{\alpha}_l^s}$ and $q_{\hat{\alpha}_u^s}$ are determined via the single parametric bootstrap procedure described above.
\label{theom_pivot}
\end{theorem}
The proof of Theorem \ref{theom_pivot} is given in the online Appendix A. A special example of Theorem \ref{theom_pivot} is that when $\mathcal{G}(\bm{x_i^\prime\beta}, A, \bm\phi)$ is a normal distribution, $\mathcal{N}(\bm{x_i^\prime \beta}, A)$, using Theorem 3.2 of \cite{chatterjee2008parametric}, we have
 \begin{equation}
\mathrm{P}(\theta_i \in \hat{I}_i^s) 
= 1 - \alpha + O(m^{-3/2}).
\end{equation}

\begin{remark}
    We illustrate the existence of a pivot for the $H_i(\bm{\beta}, A)$ in a setting where the $\theta_i$ at level 2 follows a non-normal distribution. Consider the following hierarchical model, which generalizes the model  in \citet{ghosh1987robust} by incorporating covariates in the mean structure:
    \begin{eqnarray*}
         &\text{Level 1}:&\; y_i|\theta_i, h \overset{\mathrm{ind}}{\sim} \mathcal{N}(\theta_i, h^{-1});\\
        &\text{Level 2}:& \;\theta_i|h \overset{\mathrm{ind}}{\sim} \mathcal{N}(\bm{x}_i^\prime \bm{\beta}, (\eta h)^{-1}),
    \end{eqnarray*}
    where $h$ has a proper prior distribution. Assume $h \sim \Gamma (\alpha, \lambda)$ with known shape parameter $\alpha>1$ and known scale parameter $\lambda$. In this case, the marginal distribution of $\theta_i$ is a $t$ distribution with $2\alpha$ degrees of freedom, mean $\bm{x}_i^\prime\bm{\beta}$, and scale parameter $\sqrt{\lambda/(\eta\alpha)}$: $\theta_i \sim t_{2\alpha}(\bm{x}_i^\prime\bm{\beta}, \sqrt{\lambda/(\eta\alpha)})$. Moreover, the $H_i(\bm{\beta}, A)$ also follows a $t$ distribution with known parameters:
   \begin{eqnarray*}
    H_i(\beta, A) \sim t_{2\alpha}(0, \lambda\sqrt{(\alpha-1)/\alpha}).
   \end{eqnarray*}
   Since the distribution of $H_i(\bm{\beta}, A)$ does not depend on unknown parameters, it is a pivot.
   \end{remark}

\begin{Proposition}
Given a non-pivot $H_i(\bm\beta,A)$, 
and under regularity conditions R1-R6,
we have 
    \begin{equation}
\mathrm{P}(\theta_i \in \hat{I}_i^s) 
= 1 - \alpha + O(m^{-1}),
\label{np_cll}
\end{equation}
where $\hat{I}_i^s$ is obtained from the single parametric bootstrap procedure described above. 
\label{Proposition1}
\end{Proposition}
\noindent The proof is given in the online Appendix B.

\begin{Proposition}
Suppose that:
\begin{itemize}
    \item[(i)] The random effects $u_i$ are symmetrically distributed;
    \item[(ii)] $\partial F_{2i}(x; \bm{\nu}) / \partial \nu_i >0, \; 1 \leq i \leq s$ \text{and} $\lambda_{\rm{max}}(\partial^2 F_{2i}(x; \bm{\nu})/ \partial\bm{\nu} \partial \bm{\nu}^\prime) < 0$, where $\lambda_{\rm{max}}$ means the largest eigenvalue.\\
    This condition is satisfied for some continuous distributions. For instance, when $u_i$ follows a logistic or $t$ distribution with known degrees of freedom. The only unknown parameter of $F_{2i}(x; \bm{\nu})$ is the variance $A$. As indicated in Remark \ref{moment4}, the kurtosis of $H_i(\bm{\beta}, A)$ is a decreasing function of $A$, and thus it is not difficult to show that $\partial F_{2i}(x; A) /\partial A > 0$ and $\partial^2 F_{2i}(x; A) /\partial A^2 < 0$;
    \item[(iii)] The estimators of the unknown parameters, $\hat{\bm{\nu}} = (\hat{\nu}_1, \cdots, \hat{\nu}_s)$, are either second-order unbiased or negatively biased, that is 
    $\text{E}(\hat{\nu}_i - \nu_i) = b_i + o(m^{-1})$ with $b_i \leq 0$, for $i = 1, \cdots, s$.
\end{itemize}
Under these conditions 
with the regularity conditions R1-R6
, the prediction interval \eqref{np_cll} exhibits an overcoverage property. More specifically, we can rewrite \eqref{np_cll}  as below:
    \begin{equation}
\mathrm{P}(\theta_i \in \hat{I}_i^s) 
= 1 - \alpha + T_1 + o(m^{-1}),
\label{sign_np_cll}
\end{equation}
where $T_1$ is of the order $O(m^{-1})$ with $T_1 > 0$. 

\label{overcov}
\end{Proposition}
\noindent See detailed proof in the online Appendix C. The proposition indicates that under the regularity conditions, the prediction intervals conducted by the proposed single parametric bootstrap can produce higher coverage than the nominal coverage with a non pivot $H_i(\bm{\beta}, A)$ up to the order of $O(m^{-1})$, which could be beneficial for practitioners without considering other properties of prediction intervals.

\begin{remark}
\label{moment4}
When $\mathcal{G}(\bm{x_i^\prime \beta}, A, \bm{\phi})$ is not a normal distribution, it is challenging to obtain the explicit form of the distribution of $H_i(\bm{\beta}, A)$. Consequently, it is difficult to verify whether $H_i(\bm\beta,A)$ is a pivot. Note the fact that for $H_i(\bm{\beta}, A)$ to be a pivot, its moments should not depend on any unknown parameters. Based on that, we develop a simple moment-based method to claim non-existence of pivot. Under the symmetry assumption of $u_i$, the odd moments of $H_i(\bm\beta,A)$ are zero if they exist, and the second moment equals to 1 because it is standardized. To verify if $H_i(\bm{\beta}, A)$ is a pivot, we calculate its fourth moment as follows:
\begin{eqnarray*}
    \text{E}\left[(\theta_i - \tilde{\theta}_i)/\sqrt{g_{1i}}\right]^4 
    &=& \frac{1}{g_{1i}^2}\text{E}\left[\bm{x}_i^\prime\bm{\beta} + u_i - (1-B_i)y_i - B_i\bm{x}_i^\prime\bm{\beta}\right]^4\\
    &=& 3 + \frac{D_i^2}{(A+D_i)^2}\left[\frac{\text{E}(u_i^4)}{A^2} - 3\right],
    \label{4moment_t}
\end{eqnarray*}
\noindent where $\left[\text{E}(u_i^4)/A^2 - 3\right]$ is the excess kurtosis of $u_i$. Note that when $u_i$ is  normally distributed, $\left[\text{E}(u_i^4)/A^2 - 3\right]$ is zero. When the distribution of $u_i$  is other than normal, such as t, double exponential, and logistic, $\left[\text{E}(u_i^4)/A^2 - 3\right]$ is a constant other than zero, indicating that the fourth moment of $H_i(\bm\beta,A)$ depends on the unknown parameter $A$ and thus $H_i(\bm\beta,A)$ is not a pivot. Moreover, in these cases, the fourth moment of $H_i(\bm\beta,A)$ is a decreasing function of $A$. 
\end{remark}

\section{Double parametric bootstrap}\label{sec:4}

\cite{hall2006parametric} considered parametric bootstrap methods to approximate the distribution of $(\theta_i - \bm{x}_i^\prime \hat{\bm{\beta}})/\sqrt{\hat{A}}$. Then, the prediction interval can be constructed as $(\bm{x}_i^\prime\hat{\bm\beta}+b_{\hat{\alpha}_l^s}\sqrt{\hat{A}}, \bm{x}_i^\prime\hat{\bm\beta} + b_{\hat{\alpha}_u^s}\sqrt{\hat{A}})$, where $b_{\hat{\alpha}_l^s}$ and $b_{\hat{\alpha}_u^s}$ are obtained from single bootstrap approximation based on their algorithm. Their prediction interval is based on a synthetic model or the regression model, which does not permit approximation of the conditional distribution of $\theta_i$ given the data $\bm{Y}$. As a consequence, it is likely to underweight the area specific data. 
When the level 2 distribution $\mathcal{G}$ determined only by $\bm{x}_i^\prime \bm{\beta}$ and $A$, it is easy to know that the quantity, $(\theta_i - \bm{x}_i^\prime \bm{\beta})/\sqrt{A}$, is a pivot. As \cite{hall2006parametric} stated, their prediction interval is effective as $\mathrm{P}(\bm{x}_i^\prime\hat{\bm\beta}+b_{\hat{\alpha}_l^s}\sqrt{\hat{A}}\leq \theta_i \leq \bm{x}_i^\prime\hat{\bm\beta} + b_{\hat{\alpha}_u^s}\sqrt{\hat{A}}) = 1-\alpha + O(m^{-2})$. However, when additional parameters involved into the distribution of random effects, $(\theta_i - \bm{x}_i^\prime \bm{\beta})/\sqrt{A}$ might not a pivot and we may lose the effectiveness. 

\begin{Proposition}
When considering a general distribution $\mathcal{G}(\bm{x}_i^\prime\bm\beta, A, \bm\phi)$, such as a t distribution with unknown degree of freedom, we have 
\begin{equation}
    \mathrm{P}(\bm{x}_i^\prime\hat{\bm\beta}+b_{\hat{\alpha}_l^s}\sqrt{\hat{A}}\leq \theta_i \leq \bm{x}_i^\prime\hat{\bm\beta}+b_{\hat{\alpha}_u^s}\sqrt{\hat{A}}) = 1-\alpha + O(m^{-1}),
    \label{hm_remark}
\end{equation}
where $b_{\hat{\alpha}_l^s}$ and $b_{\hat{\alpha}_u^s}$ are obtained from single bootstrap approximation based on \cite{hall2006parametric} algorithm. 
    \label{remark1}
\end{Proposition}
\noindent See proof in the online Appendix D, which is similar to Proposition \ref{Proposition1}. 

\cite{hall2006parametric} proposed a double-bootstrap method to calibrate $(\hat{\alpha}_l^s, \hat{\alpha}_u^s)$ from single parametric bootstraps, which can achieve a high degree of coverage accuracy $O(m^{-3})$. However, their calibration approach can overcorrect and produce a calibrated $\hat{\alpha}$ greater than 1, which makes the implementation infeasible. 

While the order $O(m^{-1})$ term $T_1$ in \eqref{sign_np_cll} is theoretically positive under certain conditions, it remains unclear whether this positiveness holds when  $u_i$ is asymmetrically distributed. 
Unlike \cite{hall2006parametric}, we introduce a new double bootstrap method, which does not require a pivot or symmetric $u_i$. In general, this method reduces the coverage error to $o(m^{-1})$, even when $u_i$ is asymmetrically distributed.
Our double bootstrap approach is based on the algorithm from \cite{shi1992accurate}, where the double bootstrap is proposed to obtain accurate and efficient confidence interval for parameters of interest in both nonparametric and univariate parametric distribution settings. Later on, \cite{mccullough1998implementing} discussed the theory of the double bootstrap both with and without pivoting, and provided the implementations in some nonlinear production functions. In this paper, we develop the double bootstrap in the context of our mixed effect model and apply it to obtain the EBL prediction intervals of $\theta_i$. The framework of our double parametric bootstrap is as below:
\begin{itemize}
    \item[1.] \textit{\textbf{First-stage bootstrap}}
    \begin{itemize}
    \item[1.1] Conditionally on the data $(\bm{X, Y})$, draw $\theta_i^*, \, i=1, \cdots, m$, from the distribution $\mathcal{G}(\bm{x_{i}}^\prime\hat{\bm \beta}, \hat{A}, \hat{\bm\phi})$; 
    \item[1.2] Given $\theta_i^*$, draw $y_i^*$ from the distribution $\mathcal{N}(\theta_i^*, D_i)$;
    \item[1.3] Compute $\hat{\bm \beta}^*$, $\hat{A}^*$, $\hat{\bm\phi}^*$ from the data $(\bm{X}, \bm{Y}^*)$, and obtain $\hat{\theta}_i^*$ and $\hat{g}_{1i}^*$;
     \end{itemize}
    \item[2.] \textit{\textbf{Second-stage bootstrap}}
   \begin{itemize}
    \item[2.1] Given $\hat{\bm \beta}^*$, $\hat{A}^*$, $\hat{\bm\phi}^*$, draw $\theta_i^{**}$ from the distribution $\mathcal{G}(\bm{x_{i}}^\prime\hat{\bm{\beta}}^*, \hat{A}^*, \hat{\bm\phi}^*)$; 
    \item[2.2] Given $\theta_i^{**}$, draw $y_i^{**}$ from the distribution $\mathcal{N}(\theta_i^{**}, D_i)$;
    \item[2.3] Compute $\hat{\bm \beta}^{**}$, $\hat{A}^{**}$ from the data $(\bm{X}, \bm{Y}^{**})$, also obtain $\hat{\theta}_i^{**}$ and $\hat{g}_{1i}^{**}$;
    \item[2.4] \label{2.4} Consider to obtain a $(1-\alpha)$-level, two-sided, equal-tailed prediction interval. Define
    \begin{equation}
    \hat{G}^*(z) \equiv \mathrm{P}^{**}(H_i(\hat{\bm\beta}^{**}, \hat{A}^{**}) \leq z \,| \, \hat{\bm \beta}^*, \hat{A}^*, \hat{\phi}^*)
    \label{G:def}
    \end{equation}
    For seeking the upper limit, we firstly solve the following system of equations in order to obtain $\hat{\alpha}_u$ such that
    \begin{equation}
  \left\{\begin{array}{@{}l@{}}
    \mathrm{P}^*(H_i(\hat{\bm\beta}^*, \hat{A}^*) \leq q_{\hat{\alpha}_u} \,| \, \hat{\bm \beta}, \hat{A}, \hat{\bm\phi}) = 1 - \alpha/2\\
    \mathrm{P}^{**}(H_i(\hat{\bm\beta}^{**}, \hat{A}^{**}) \leq q_{\hat{\alpha}_u} \,| \, \hat{\bm \beta}^*, \hat{A}^*, \hat{\bm\phi}^*) = {\hat{\alpha}_u}.
  \end{array}\right.\,,
  \label{db}
\end{equation}
Using the definition \eqref{G:def}, we have $q_{\hat{\alpha}_u} = \hat{G}^{*-1}(\hat{\alpha}_u)$. Then, the above system of equations is equivalent to
   \begin{equation}
   \mathrm{P}^*\left(H_i(\hat{\bm\beta}^*, \hat{A}^*) \leq \hat{G}^{*-1}(\hat{\alpha}_u) \,| \, \hat{\bm \beta}, \hat{A}, \hat{\bm\phi}\right) = 1 - \alpha/2.
   \label{eqv:sys}
   \end{equation}
Rewriting \eqref{eqv:sys} gives
    \begin{equation}
\mathrm{P}^*\left[\mathrm{P}^{**}\left(H_i(\hat{\bm\beta}^{**}, \hat{A}^{**}) \leq H_i(\hat{\bm\beta}^*, \hat{A}^*) \,| \, \hat{\bm \beta}^*, \hat{A}^*, \hat{\bm\phi}^*\right) \leq \hat{\alpha}_u \,| \, \hat{\bm \beta}, \hat{A}, \hat{\bm\phi}\right] = 1 - \alpha/2.
   \end{equation}
   Note that the inner probability, $\mathrm{P}^{**}(H_i(\hat{\bm\beta}^{**}, \hat{A}^{**}) \leq H_i(\hat{\bm\beta}^*, \hat{A}^*) \,| \, \hat{\bm \beta}^*, \hat{A}^*, \hat{\bm\phi}^*)$, is a function of the first-stage bootstrap resample. More specifically, on the $j$th first-stage bootstrap sample, after all $K$ second-stage bootstrap operations are completed, let $Z_j$ be the proportion of times that $H_i(\hat{\bm\beta}^{**}, \hat{A}^{**}) \leq H_i(\hat{\bm\beta}^*, \hat{A}^*)$, i.e.,
       \begin{eqnarray*}
         Z_j = \mathrm{P}^{**}\left(H_i(\hat{\bm\beta}^{**}, \hat{A}^{**}) \leq H_i(\hat{\bm\beta}^*, \hat{A}^*) \,| \, \hat{\bm \beta}^*, \hat{A}^*, \hat{\bm\phi}^*\right) \approx \#(H_i(\hat{\bm\beta}^{**}, \hat{A}^{**}) \leq H_i(\hat{\bm\beta}^*, \hat{A}^*)) / K
       \end{eqnarray*}
We will use $Z_j$ to adjust the first-stage intervals. After all bootstrapping operations are complete, we have estimates $Z_j$, $j = 1, 2, \cdots, J$, where $J$ is the number of first-stage bootstrap operations. Sort $Z_j$ and choose the $1-\alpha/2$ quantile as the percentile point $\hat{\alpha}_u$ for defining the double-bootstrap upper limit.
 
    \item[2.5] After completing the step 2.4, we have all estimates $H_i^{(j)}(\hat{\bm\beta}^*, \hat{A}^*)$, $j = 1, \cdots, J$ and $\hat{\alpha}_u$. Choose $q_{\alpha_u} = q_{\hat{\alpha}^d_u}$ such that
    \begin{equation}
    \mathrm{P}^*(H_i(\hat{\bm\beta}^*, \hat{A}^*) \leq q_{\alpha_u} \,| \, \hat{\bm \beta}, \hat{A}, \hat{\bm\phi}) = \hat{\alpha}_u.
    \label{sbdb}
    \end{equation}
    \item[2.6] Take $q_{\hat{\alpha}^d_u}$ as the upper limit of the double bootstrap prediction interval. Similar operations determine the lower limit, $q_{\hat{\alpha}^d_l}$. Finally, construct the prediction interval of $\theta_i$ as $$ \hat{I}_i^d = \left(\hat{\theta}_i + q_{\hat{\alpha}^d_l}\sqrt{\hat{g}_{1i}} \, , \, \hat{\theta}_i + q_{\hat{\alpha}^d_u}\sqrt{\hat{g}_{1i}}\right).$$
     \end{itemize}    
\end{itemize}

The algorithm above shows that the single parametric bootstrap $\hat{I}^s_i$ is calibrated by the second-stage bootstrap. Such calibration improves the coverage accuracy to $o(m^{-1})$ even when the pivot does not exist. One more advantage of our double bootstrap algorithm is that it avoids the problem of over correction and so make it more practical than that of \cite{hall2006parametric}.
\begin{theorem}
Under regularity conditions, for a preassigned $\alpha \in (0, 1)$ and arbitrary $i = 1, \cdots, m$,  we have
\begin{equation}
\mathrm{P}\left(\hat{\theta}_i+q_{\hat{\alpha}^d_l}\sqrt{\hat{g}_{1i}} \leq \theta_i \leq \hat{\theta}_i+q_{\hat{\alpha}_u^d}\sqrt{\hat{g}_{1i}}\right) = 1 - \alpha + o(m^{-1}),
\end{equation}
where $q_{\hat{\alpha}^d_l}$ and $q_{\hat{\alpha}^d_u}$ are obtained from the double parametric bootstrap procedure described above.
\label{theom_nonpivot}
\end{theorem}
See proof of the theorem in online Appendix E.

\section{Monte Carlo Simulations}\label{sec:5}
In this section, we compare the performance of the proposed parametric bootstrap with their competitors where available, using Monte Carlo simulation studies. To maintain comparability with existing studies, we adopt part of the simulation framework of \cite{chatterjee2008parametric}. We consider an area-level model with $\bm{x}_i^\prime\bm{\beta} = 0$, and five groups of small areas. Within each group, the sampling variances $D_i$'s remain the same. Two patterns for the $D_i$'s are considered: (i) $(4.0, 0.6, 0.5, 0.4, 0.2)$; (ii) $(8.0, 1.2, 1.0, 0.8, 0.4)$. For pattern (i), we take $A = 1$; for pattern (ii), we take $A=2$ doubling the variances of pattern (i) while preserving the $B_i = D_i/(A+D_i)$ ratios. To examine the effect of $m$, we consider $m = 15$ and 50. With the increase of $m$,  all methods improve and get closer to one another, supporting our asymptotic theory. Since we obtained virtually identical results under these two patterns, for the full study we confined attention to the pattern (i) and the results for pattern (ii) are provided in the Supplementary Material.

The Prasad-Rao method-of-moments \citep{prasad1990estimation}, and the Fay-Herriot method of estimating the variance component $A$ are considered. For the Fay-Herriot estimator of $A$, we employ the method of scoring to solve the estimating equation, which has showed to be more stable  \citep{datta2005measuring} than the Newton-Raphson method originally used in \cite{fay1979estimates}. The estimation equation of the Fay-Herriot estimator is
\[f(A) = \sum_{i=1}^m\frac{(y_i - \bm{x}_i^\prime\tilde{\bm\beta})^2}{A+D_i} = 0. \]
Here, $f(A)$ is a decreasing function with respect to $A$, and the expectation of its first derivative $\text{E}(f^\prime(A)) = -\rm{tr}(\bm{P}) < 0$. To improve computational efficiency, we implemented a slight modification to the algorithm. Specifically, the revised Fay-Herriot algorithm begins by calculating $f(A)$ at the initial point $A_0 = 0$, as in \cite{fay1979estimates}. If $f(A_0) < 0$, we truncate the estimator to $\hat{A}_{\rm{FH}} = 0.01$. Otherwise, the iterative process continues to search for a positive solution, with $\hat{A}_{\rm{FH}} = 0.01$ also applied if no positive solution is found. This revised Fay-Herriot algorithm further enhances computational efficiency compared to the original method.

\subsection{Simulations on symmetric cases}
First, we consider the scenario that the random effects $u_i$ are symmetrically distributed. Specifically, we assume $\mathcal{G}$ is a $t$ distribution with 9 degrees of freedom. In this setting, we compare coverage probabilities and average lengths of the following seven different prediction intervals of $\theta_i$:
\begin{itemize}
    \item Two prediction intervals based on the proposed single parametric bootstrap method with two different variance estimators: $\hat{A}_{\rm{FH}}$ and $\hat{A}_{\rm{PR}}$, denoted as \textbf{SB.FH} and \textbf{SB.PR}, respectively;
    \item Two prediction intervals based on the single parametric bootstrap methods proposed by \cite{hall2006parametric}, using the same variance estimators, denoted as \textbf{HM.FH} and \textbf{HM.PR};
    \item Two traditional prediction intervals of the form $\hat\theta_i \pm z_{\alpha/2}\sqrt{\text{mspe}(\hat{\theta}_i)}$, where the variance estimator $\hat{A}_{\rm{FH}}$ is used along with the MSPE estimator of \cite{datta2005measuring}, denoted as \textbf{FH}, and the variance estimator $\hat{A}_{\rm{PR}}$ is used with the Prasad-Rao MSPE estimator \citep{prasad1990estimation}, denoted as \textbf{PR};
    \item The direct confidence interval (\textbf{DIRECT}) given by $y_i \pm z_{\alpha/2}\sqrt{D_i}$.
\end{itemize}

Each reported result is based on 1000 simulation runs. For all cases, we consider single bootstrap sample of size 400 and three different nominal coverages: $\alpha \in (80\%, 90\%, 95\%)$. 
All computations are carried out  in \texttt{R} on a \texttt{13th Gen Intel(R) Core(TM) i7-13700F 2.10 GHz} processor. The average computation time per simulation run for the proposed single bootstrap method is 10.05 seconds for $m = 15$ and 11.55 seconds for $m = 50$.

We report the percentages of negative estimates for A under $t_9$ distribution in Table \ref{tab:z0}. For $m =15$, the Prasad-Rao method-of-moment approach yields as high as about 13\% negative estimates of $A$ under pattern (i) of $D_i$. Fay-Herroit method produces significantly fewer negative estimates than the Prasad-Rao method-of-moment approach across all scenarios.
\begin{table}[htbp]
    \centering
    \caption{Percentages of negative estimates of $A$ ($\hat{A}^*$) for different estimation methods and different patterns of $D_i$ under $t_9$.}
    \begin{tabular}{cccccc}
    \toprule
              &   \multicolumn{2}{c}{Pattern (i)} & \multicolumn{2}{c}{Pattern (ii)}\\
    \midrule
    $m$     &  FH & PR &  FH & PR\\
    15     &  1.6 (6.4) & 13.1 (25.0) & 1.1 (5.7) & 11.8 ( 23.2)\\
    50 &   0 (0.07) & 1.2(6.6) & 0 (0.07) & 1.2 (6.1)\\
    \bottomrule
    \end{tabular}
    \label{tab:z0}
\end{table}

Table \ref{tab1:t50} presents coverage probabilities and average lengths for each prediction interval method with $m = 50$ and the $t_9$ distribution under the pattern (i) of $D_i$. The \textbf{SB.PR} and \textbf{HM.PR} prediction intervals consistently over-cover. \textbf{SB.PR} prediction intervals have shorter lengths than \textbf{HM.PR}, especially for groups 2-5 with smaller sampling variances $D_i$. The \textbf{SB.FH} and \textbf{HM.FH} prediction intervals perform well regarding coverage errors. Specifically, their coverage probabilities are very close to all three nominal coverages. The \textbf{SB.FH} method produces the shortest average lengths among these four single parametric bootstrap methods. The \textbf{FH} intervals tend to undercover for G1 group when $\alpha = 90\%$ and $95\%$. The \textbf{PR} intervals have the undercoverage issue for G1 group across all three nominal coverages.

\begin{table}[!t]
\small
\centering
\caption{Coverage probabilities (average length) of different intervals under $t_9$ distribution with $m$ = 50 and Pattern (i)}
\begin{adjustbox}{width=\textwidth}
\begin{tabular}{cccccccccc}
  \toprule
 & SB.FH & HM.FH & SB.PR & HM.PR & FH & PR & DIRECT \\
  \midrule
\multicolumn{8}{l}{(i) $\alpha = 80\%$}\\
G1 & 79.94 ( 2.31 )  & 80.30 ( 2.53 )  & 83.70 ( 2.59 )  & 84.45 ( 2.84 )  & 79.89 ( 2.30 )  & 78.66 ( 2.31 )  & 79.43 ( 5.13 )  \\ 
  G2 & 80.15 ( 1.59 )  & 80.31 ( 2.49 )  & 83.63 ( 1.79 )  & 84.54 ( 2.81 )  & 80.38 ( 1.59 )  & 81.55 ( 1.67 )  & 79.85 ( 1.99 )  \\ 
  G3 & 79.71 ( 1.50 )  & 79.53 ( 2.48 )  & 83.12 ( 1.69 )  & 83.75 ( 2.80 )  & 79.86 ( 1.50 )  & 80.88 ( 1.58 )  & 79.99 ( 1.81 )  \\ 
  G4 & 80.47 ( 1.39 )  & 80.03 ( 2.48 )  & 84.06 ( 1.57 )  & 83.95 ( 2.80 )  & 80.49 ( 1.39 )  & 82.06 ( 1.47 )  & 80.29 ( 1.62 )  \\ 
  G5 & 79.65 ( 1.05 )  & 80.48 ( 2.46 )  & 83.32 ( 1.20 )  & 84.89 ( 2.78 )  & 79.99 ( 1.06 )  & 81.8 ( 1.14 )  & 80.32 ( 1.15 )  \\ 
  \multicolumn{8}{l}{(i) $\alpha = 90\%$}\\
  G1 & 89.88 ( 3.05 )  & 90.35 ( 3.39 )  & 94.13 ( 3.79 )  & 94.25 ( 4.21 )  & 88.90 ( 2.96 )  & 87.82 ( 2.97 )  & 89.42 ( 6.58 )  \\ 
  G2 & 90.17 ( 2.05 )  & 90.37 ( 3.30 )  & 93.25 ( 2.59 )  & 94.27 ( 4.14 )  & 89.89 ( 2.04 )  & 90.58 ( 2.14 )  & 90.19 ( 2.55 )  \\ 
  G3 & 90.06 ( 1.94 )  & 89.34 ( 3.30 )  & 93.38 ( 2.44 )  & 93.76 ( 4.14 )  & 90.15 ( 1.92 )  & 90.78 ( 2.03 )  & 90.11 ( 2.33 )  \\ 
  G4 & 90.07 ( 1.79 )  & 89.50 ( 3.29 )  & 93.41 ( 2.26 )  & 94.24 ( 4.13 )  & 90.12 ( 1.78 )  & 91.19 ( 1.89 )  & 90.24 ( 2.08 )  \\ 
  G5 & 90.29 ( 1.36 )  & 90.28 ( 3.26 )  & 93.33 ( 1.72 )  & 94.2 ( 4.09 )  & 90.26 ( 1.36 )  & 91.56 ( 1.47 )  & 90.32 ( 1.47 )  \\ 
  \multicolumn{8}{l}{(i) $\alpha = 95\%$}\\
  G1 & 95.13 ( 3.75 )  & 95.04 ( 4.22 )  & 97.30 ( 5.31 )  & 97.42 ( 5.91 )  & 93.57 ( 3.52 )  & 92.65 ( 3.53 )  & 94.77 ( 7.84 )  \\ 
  G2 & 95.11 ( 2.47 )  & 95.35 ( 4.07 )  & 96.99 ( 3.61 )  & 97.5 ( 5.79 )  & 94.87 ( 2.43 )  & 95.52 ( 2.55 )  & 95.38 ( 3.04 )  \\ 
  G3 & 94.96 ( 2.32 )  & 94.45 ( 4.06 )  & 97.05 ( 3.39 )  & 97.42 ( 5.77 )  & 94.70 ( 2.29 )  & 95.51 ( 2.41 )  & 95.00 ( 2.77 )  \\ 
  G4 & 95.06 ( 2.14 )  & 95.03 ( 4.04 )  & 96.86 ( 3.14 )  & 97.58 ( 5.75 )  & 95.12 ( 2.12 )  & 95.86 ( 2.25 )  & 95.07 ( 2.48 )  \\ 
  G5 & 95.31 ( 1.62 )  & 94.98 ( 3.99 )  & 96.87 ( 2.38 )  & 97.42 ( 5.69 )  & 95.28 ( 1.62 )  & 96.01 ( 1.75 )  & 95.17 ( 1.75 )  \\ 
   \bottomrule
\end{tabular}
\label{tab1:t50}
\end{adjustbox}
\end{table}

Table \ref{tab3:t15} reports the results for $m = 15$. As illustrated in Table \ref{tab:z0}, the Prasad–Rao method produces an extremely high percentage of zero estimates for a small $m = 15$. Thus, it is not surprising to note that the \textbf{SB.PR} intervals have severe undercoverage problem when $\alpha = 90\%$ and 95\%. The high percentage of negative estimates of $A$ might also contribute to the similar undercoverage problem of \textbf{HM.PR} intervals at $\alpha = 90\%$ and 95\% as well as the large average lengths of both \textbf{SB.PR} and \textbf{HM.PR} intervals. The \textbf{SB.FH} and \textbf{HM.FH} intervals uniformly tend to overcover. Still, \textbf{SB.FH} intervals have the shortest average lengths among the four types of parametric bootstrap intervals. The \textbf{FH} intervals significantly undercover for group G1 at three normial coverages and have slight undercoverage for groups G2 and G3, when $\alpha = 90\%$ and 95\%. The \textbf{PR} intervals undercover for group G1 and switch to overcover for rest of the groups at all three nominal coverages. 

\begin{table}[!t]
\small
\centering
\caption{Coverage probabilities (average length) of different intervals under $t_9$ distribution with $m$ = 15 and Pattern (i)}
\begin{adjustbox}{width=\textwidth}
\begin{tabular}{ccccccccc}
  \toprule
 & SB.FH & HM.FH & SB.PR & HM.PR & FH & PR & DIRECT \\
  \midrule
  \multicolumn{8}{l}{(i) $\alpha = 80\%$}\\
G1 & 82.43 ( 2.68 )  & 82.50 ( 2.94 )  & 81.50 ( 3.37 )  & 82.30 ( 3.74 )  & 75.07 ( 2.31 )  & 76.23 ( 2.45 )  & 81.37 ( 5.13 )  \\ 
  G2 & 82.80 ( 1.78 )  & 83.50 ( 2.78 )  & 80.77 ( 2.25 )  & 82.13 ( 3.60 )  & 79.43 ( 1.63 )  & 88.97 ( 2.19 )  & 80.03 ( 1.99 )  \\ 
  G3 & 81.40 ( 1.67 )  & 81.37 ( 2.76 )  & 80.57 ( 2.12 )  & 82.27 ( 3.58 )  & 79.07 ( 1.54 )  & 88.20 ( 2.15 )  & 80.40 ( 1.81 )  \\ 
  G4 & 84.07 ( 1.55 )  & 84.17 ( 2.74 )  & 81.50 ( 1.95 )  & 82.97 ( 3.54 )  & 81.40 ( 1.43 )  & 90.47 ( 2.1 )  & 80.67 ( 1.62 )  \\ 
  G5 & 81.80 ( 1.17 )  & 83.17 ( 2.66 )  & 79.33 ( 1.47 )  & 82.23 ( 3.42 )  & 80.73 ( 1.11 )  & 88.50 ( 1.98 )  & 80.00 ( 1.15 )  \\ 
  \multicolumn{8}{l}{(i) $\alpha = 90\%$}\\
  G1 & 93.07 ( 3.89 )  & 92.87 ( 4.33 )  & 87.00 ( 5.77 )  & 87.3 ( 6.45 )  & 83.97 ( 2.97 )  & 85.80 ( 3.14 )  & 90.33 ( 6.58 )  \\ 
  G2 & 93.00 ( 2.51 )  & 93.1 ( 3.96 )  & 86.30 ( 3.78 )  & 87.27 ( 6.03 )  & 89.37 ( 2.09 )  & 95.70 ( 2.81 )  & 90.87 ( 2.55 )  \\ 
  G3 & 92.63 ( 2.35 )  & 92.5 ( 3.91 )  & 86.50 ( 3.57 )  & 87.77 ( 6.01 )  & 88.37 ( 1.98 )  & 95.27 ( 2.76 )  & 90.20 ( 2.33 )  \\ 
  G4 & 93.90 ( 2.16 )  & 93.83 ( 3.86 )  & 87.53 ( 3.28 )  & 87.87 ( 5.94 )  & 90.17 ( 1.84 )  & 96.67 ( 2.70 )  & 90.93 ( 2.08 )  \\ 
  G5 & 91.83 ( 1.60 )  & 92.87 ( 3.67 )  & 86.03 ( 2.42 )  & 87.53 ( 5.62 )  & 90.03 ( 1.43 )  & 94.73 ( 2.54 )  & 89.60 ( 1.47 )  \\
  \multicolumn{8}{l}{(i) $\alpha = 95\%$}\\
  G1 & 97.13 ( 5.46 )  & 97.07 ( 6.11 )  & 88.77 ( 8.91 )  & 88.80 ( 9.98 )  & 88.73 ( 3.54 )  & 90.73 ( 3.75 )  & 94.73 ( 7.84 )  \\ 
  G2 & 96.33 ( 3.41 )  & 96.83 ( 5.36 )  & 88.43 ( 5.73 )  & 88.87 ( 9.17 )  & 93.7 ( 2.49 )  & 97.87 ( 3.35 )  & 95.7 ( 3.04 )  \\ 
  G3 & 96.60 ( 3.18 )  & 96.83 ( 5.26 )  & 89.20 ( 5.39 )  & 89.50 ( 9.10 )  & 93.33 ( 2.36 )  & 97.90 ( 3.29 )  & 94.77 ( 2.77 )  \\ 
  G4 & 97.23 ( 2.90 )  & 97.37 ( 5.15 )  & 89.87 ( 4.93 )  & 89.90 ( 8.95 )  & 95.33 ( 2.19 )  & 98.83 ( 3.22 )  & 95.63 ( 2.48 )  \\ 
  G5 & 96.30 ( 2.11 )  & 96.93 ( 4.77 )  & 88.93 ( 3.59 )  & 89.30 ( 8.38 )  & 95.13 ( 1.7 )  & 97.80 ( 3.03 )  & 94.57 ( 1.75 )  \\ 
   \bottomrule
\end{tabular}
\label{tab3:t15}
\end{adjustbox}
\end{table}

\subsection{Further simulations on asymmetric cases}
While some of our theoretical results for single bootstrap are based on the symmetry assumption, we also use simulation studies to assess our proposed parametric bootstrap method under asymmetric conditions. Specifically, for the same models and parameter choices as above, we change $\mathcal{G}$ to be a shifted exponential distribution (SE). Besides the seven prediction intervals mentioned before, we also include our proposed double parametric bootstrap intervals in this subsection and we expect they might improve the potential coverage error under asymmetric conditions. Below, the double bootstrap method based on the Fay-Herriot variance estimator is denoted by \textbf{DB.FH}, and the double bootstrap method based on the Prasad-Rao variance estimator is indicated by \textbf{DB.PR}. We keep bootstrap sample of size 400 in the first stage and apply two different sizes $B_2 \in (50, 100)$ in the second stage. Since for $B_2 =50$ and $B_2 =100$, we obtained very similar results. Moreover, the two patterns of $D_i$ gave us similar results. In the following, we provide detailed discussions only for $B_2 =100$ under the pattern (i), but the results for $B_2 =50$ are provided in the Supplementary Material. 
The average computation time per simulation run for the proposed double bootstrap method is 22.41 seconds for $m = 15,\; B_2 = 50$; 36.61 seconds for $m = 15,\; B_2 = 100$; 57.60 seconds for $m = 50,\; B_2 = 50$; and 97.10 seconds for $m = 50,\; B_2 = 100$.

Table \ref{tab2:z0} shows the percentages of negative estimates for $A$ under the SE distribution. When $m = 15$, the Prasad-Rao method-of-moments approach yields very high percentages of negative estimates at both the first and second stages of bootstrapping. Although the Fay-Herriot method results in fewer negative estimates than the Prasad-Rao approach, the occurrence of negative estimates remains notable, especially when the number of small areas is small (e.g., $m = 15$).
\begin{table}[htbp]
    \centering
    \caption{Percentages of negative estimates of $A$ ($\hat{A}^*$)[$\hat{A}^{**}$] for different estimation methods and different patterns of $D_i$ under SE distribution.}
    \begin{tabular}{cccccc}
    \toprule
              &   \multicolumn{2}{c}{Pattern (i)} & \multicolumn{2}{c}{Pattern (ii)}\\
              \midrule
    $m$     &  FH & PR &  FH & PR\\
    15     &  3.6 (10.72) [16.32] &  16.7 (27.98) [31.95] & 3.9 (10.93) [16.46] & 18.1 (28.49) [31.97]\\
    50 &   0 (0.37) [1.44] &  1.80 ( 7.89) [ 12.82] & 0 (0.44) [1.56] & 2.80 (8.89) [13.68]\\
    \bottomrule
    \end{tabular}
    \label{tab2:z0}
\end{table}

Table \ref{tab1:se50} displays the coverage probabilities and average lengths for each prediction interval under $m = 50$. The parametric bootstrap methods based on Prasad-Rao estimators of $A$, \textbf{SB.PR}, \textbf{HM.PR} and \textbf{DB.PR}, consistently tend to over-cover at $\alpha = 80\%$ and 90\%. \textbf{DB.PR} intervals bring the coverage probabilities closer to the nominal coverages but have larger lengths than single bootstrap \textbf{SB.PR} intervals. Even with $m=50$, the relatively high percentage of zero estimates of $A$ in the second bootstrap stage showed in Table \ref{tab2:z0}, might contribute to the increment on interval length of \textbf{DB.PR}. 
When applying the Fay-Herriot method in estimating $A$, our proposed single parametric bootstrap intervals \textbf{SB.FH} already showed good performance, \textbf{DB.FH} has little or no effect. The \textbf{FH} intervals perform well, except overcoverage at $\alpha = 80\%$. \textbf{PR} intervals have undercoverage problem for group G1 at $\alpha = 95\%$. Overall, the \textbf{SB.FH} and \textbf{FH} intervals perform the best in terms of both coverage probabilities and average lengths in this setting.

\begin{table}[!t]
\small
\centering
\caption{Coverage probabilities (average length) of different intervals under shifted exponential distribution with $m$ = 50, $B_2$ = 100 and Pattern (i)}
\begin{adjustbox}{width=\textwidth}
\begin{tabular}{cccccccccccc}
  \toprule
 & SB.FH & HM.FH &DB.FH & SB.PR & HM.PR &DB.PR & FH & PR & DIRECT \\
  \midrule
  \multicolumn{10}{l}{(i) $\alpha = 80\%$}\\
G1 & 80.29 ( 2.20 )  & 80.75 ( 2.35 )  & 79.58 ( 2.22 )  & 84.50 ( 2.50 )  & 87.12 ( 2.68 )  & 83.16 ( 2.86 )  & 83.45 ( 2.28 )  & 81.94 ( 2.30 )  & 79.90 ( 5.13 )  \\ 
  G2 & 79.99 ( 1.55 )  & 81.01 ( 2.31 )  & 79.42 ( 1.56 )  & 83.66 ( 1.77 )  & 87.17 ( 2.65 )  & 82.99 ( 2.07 )  & 80.56 ( 1.57 )  & 82.99 ( 1.67 )  & 79.94 ( 1.99 )  \\ 
  G3 & 80.06 ( 1.47 )  & 81.27 ( 2.30 )  & 79.62 ( 1.47 )  & 83.41 ( 1.68 )  & 87.39 ( 2.65 )  & 83.09 ( 1.97 )  & 80.38 ( 1.48 )  & 82.61 ( 1.59 )  & 79.86 ( 1.81 )  \\ 
  G4 & 80.60 ( 1.36 )  & 80.68 ( 2.30 )  & 80.20 ( 1.37 )  & 83.71 ( 1.56 )  & 87.00 ( 2.64 )  & 82.96 ( 1.83 )  & 81.05 ( 1.37 )  & 83.78 ( 1.48 )  & 80.15 ( 1.62 )  \\ 
  G5 & 80.81 ( 1.05 )  & 80.52 ( 2.29 )  & 80.34 ( 1.05 )  & 83.79 ( 1.20 )  & 87.13 ( 2.63 )  & 83.23 ( 1.4 )  & 80.95 ( 1.05 )  & 83.54 ( 1.17 )  & 80.69 ( 1.15 )  \\ 
  \multicolumn{10}{l}{(i) $\alpha = 90\%$}\\
 G1 & 90.46 ( 3.01 )  & 91.81 ( 3.26 )  & 89.56 ( 3.03 )  & 93.49 ( 3.79 )  & 94.70 ( 4.16 )  & 92.15 ( 5.53 )  & 90.98 ( 2.92 )  & 90.23 ( 2.95 )  & 90.01 ( 6.58 )  \\ 
  G2 & 90.22 ( 2.04 )  & 91.74 ( 3.16 )  & 89.53 ( 2.05 )  & 93.01 ( 2.61 )  & 94.91 ( 4.10 )  & 92.09 ( 4.05 )  & 90.14 ( 2.02 )  & 92.22 ( 2.14 )  & 89.81 ( 2.55 )  \\ 
  G3 & 90.00 ( 1.93 )  & 91.93 ( 3.15 )  & 89.48 ( 1.93 )  & 92.51 ( 2.46 )  & 94.89 ( 4.10 )  & 91.69 ( 3.85 )  & 90.2 ( 1.9 )  & 91.69 ( 2.03 )  & 89.80 ( 2.33 )  \\ 
  G4 & 90.48 ( 1.78 )  & 91.74 ( 3.14 )  & 89.8 ( 1.79 )  & 92.72 ( 2.28 )  & 94.65 ( 4.09 )  & 92.09 ( 3.59 )  & 90.47 ( 1.76 )  & 92.06 ( 1.90 )  & 90.17 ( 2.08 )  \\ 
  G5 & 90.69 ( 1.36 )  & 91.34 ( 3.10 )  & 90.24 ( 1.37 )  & 92.97 ( 1.74 )  & 94.64 ( 4.04 )  & 92.23 ( 2.76 )  & 90.71 ( 1.35 )  & 92.35 ( 1.50 )  & 90.29 ( 1.47 )  \\
  \multicolumn{10}{l}{(i) $\alpha = 95\%$}\\
G1 & 95.43 ( 3.83 )  & 96.68 ( 4.18 )  & 94.37 ( 3.90 )  & 96.45 ( 5.44 )  & 96.80 ( 6.02 )  & 95.73 ( 9.56 )  & 94.27 ( 3.48 )  & 93.58 ( 3.52 )  & 95.04 ( 7.84 )  \\ 
  G2 & 95.32 ( 2.51 )  & 96.69 ( 3.98 )  & 94.77 ( 2.58 )  & 96.37 ( 3.69 )  & 96.91 ( 5.90 )  & 95.88 ( 6.82 )  & 94.91 ( 2.40 )  & 96.04 ( 2.55 )  & 95.07 ( 3.04 )  \\ 
  G3 & 94.87 ( 2.36 )  & 96.54 ( 3.96 )  & 94.35 ( 2.43 )  & 96.34 ( 3.49 )  & 97.19 ( 5.89 )  & 95.72 ( 6.48 )  & 94.72 ( 2.27 )  & 95.51 ( 2.42 )  & 94.79 ( 2.77 )  \\ 
  G4 & 95.41 ( 2.17 )  & 96.36 ( 3.93 )  & 94.87 ( 2.23 )  & 96.35 ( 3.22 )  & 97.01 ( 5.86 )  & 95.87 ( 6.06 )  & 95.13 ( 2.10 )  & 95.92 ( 2.27 )  & 95.17 ( 2.48 )  \\ 
  G5 & 95.34 ( 1.64 )  & 96.27 ( 3.86 )  & 95.16 ( 1.69 )  & 96.34 ( 2.43 )  & 96.99 ( 5.79 )  & 95.92 ( 4.60 )  & 95.31 ( 1.61 )  & 96.30 ( 1.78 )  & 95.44 ( 1.75 )  \\
   \bottomrule
\end{tabular}
\label{tab1:se50}
\end{adjustbox}
\end{table}

With $m=15$, similar to the cases under $t_9$ distribution, the Prasad-Rao variance estimator produces very high percentages of zero estimates of $A$, especially at the bootstrap stages. See Table \ref{tab2:z0}. Thus, it is not abnormal to have the results that the three parametric bootstrap intervals based on $\hat{A}_{PR}$ have very low coverage probabilities when the nominal coverage $\alpha = 90\%$ or 95\%. See Table \ref{tab1:se15} for the cases $m=15$. The other three parametric bootstrap intervals based on $\hat{A}_{FH}$ overcover when nominal covarage is 80\% or 90\%. When $\alpha = 95\%$, \textbf{SB.FH} intervals have good performance as well as \textbf{DB.FH} in terms of coverage error, while \textbf{SB.FH} intervals have shorter average lengths than those of \textbf{HM.FH} and \textbf{DB.FH}. The \textbf{FH} intervals suffer from severe undercoverage for group G1 at all nominal coverages as well as \textbf{PR} intervals at $\alpha = 95\%$. With $\alpha = 80\%$ and 90\%, \textbf{PR} intervals have longer average lengths than those of \textbf{SB.FH} intervals for groups G2 - G5.
\begin{table}[!t]
\small
\centering
\caption{Coverage probabilities (average length) of different intervals under shifted exponential distribution with $m$ = 15, $B_2$ = 100 and Pattern (i)}
\begin{adjustbox}{width=\textwidth}
\begin{tabular}{cccccccccccc}
  \toprule
 & SB.FH & HM.FH &DB.FH & SB.PR & HM.PR &DB.PR & FH & PR & DIRECT \\
  \midrule
  \multicolumn{10}{l}{(i) $\alpha = 80\%$}\\
G1 & 84.63 ( 2.66 )  & 86.00 ( 2.90 )  & 83.63 ( 2.68 )  & 79.77 ( 3.31 )  & 81.33 ( 3.69 )  & 80.6 ( 4.89 )  & 77.17 ( 2.23 )  & 80.43 ( 2.43 )  & 79.8 ( 5.13 )  \\ 
  G2 & 82.40 ( 1.78 )  & 85.5 ( 2.74 )  & 81.63 ( 1.80 )  & 79.50 ( 2.20 )  & 81.33 ( 3.56 )  & 79.57 ( 3.26 )  & 79.43 ( 1.59 )  & 88.6 ( 2.21 )  & 79.23 ( 1.99 )  \\ 
  G3 & 82.77 ( 1.68 )  & 85.87 ( 2.72 )  & 81.97 ( 1.70 )  & 80.17 ( 2.08 )  & 81.53 ( 3.55 )  & 80.40 ( 3.07 )  & 80.13 ( 1.51 )  & 89.27 ( 2.18 )  & 78.7 ( 1.81 )  \\ 
  G4 & 84.13 ( 1.55 )  & 84.90 ( 2.70 )  & 82.97 ( 1.57 )  & 79.63 ( 1.92 )  & 80.83 ( 3.51 )  & 79.63 ( 2.80 )  & 81.83 ( 1.41 )  & 90.9 ( 2.15 )  & 82.03 ( 1.62 )  \\ 
  G5 & 82.80 ( 1.18 )  & 84.23 ( 2.61 )  & 81.70 ( 1.18 )  & 79.7 ( 1.44 )  & 80.67 ( 3.38 )  & 80.00 ( 2.00 )  & 81.40 ( 1.12 )  & 88.93 ( 2.07 )  & 80.60 ( 1.15 )  \\
  \multicolumn{10}{l}{(i) $\alpha = 90\%$}\\
G1 & 93.77 ( 4.08 )  & 94.33 ( 4.54 )  & 93.13 ( 4.83 )  & 85.83 ( 5.84 )  & 86.10 ( 6.59 )  & 86.07 ( 13.91 )  & 86.13 ( 2.86 )  & 89.37 ( 3.12 )  & 90.43 ( 6.58 )  \\ 
  G2 & 92.70 ( 2.59 )  & 94.00 ( 4.06 )  & 91.53 ( 2.93 )  & 86.00 ( 3.78 )  & 86.33 ( 6.19 )  & 85.97 ( 8.66 )  & 89.07 ( 2.04 )  & 95.13 ( 2.84 )  & 89.67 ( 2.55 )  \\ 
  G3 & 92.83 ( 2.43 )  & 93.97 ( 4.01 )  & 91.67 ( 2.69 )  & 86.80 ( 3.55 )  & 86.80 ( 6.14 )  & 86.73 ( 8.16 )  & 89.5 ( 1.94 )  & 95.77 ( 2.8 )  & 90.20 ( 2.33 )  \\ 
  G4 & 93.20 ( 2.23 )  & 93.83 ( 3.94 )  & 92.00 ( 2.48 )  & 85.93 ( 3.26 )  & 86.37 ( 6.07 )  & 85.5 ( 7.39 )  & 91.43 ( 1.81 )  & 96.2 ( 2.76 )  & 91.60 ( 2.08 )  \\ 
  G5 & 92.17 ( 1.65 )  & 93.1 ( 3.71 )  & 91.00 ( 1.75 )  & 86.53 ( 2.39 )  & 86.57 ( 5.74 )  & 86.33 ( 5.20 )  & 91.30 ( 1.43 )  & 95.23 ( 2.66 )  & 90.23 ( 1.47 )  \\ 
  \multicolumn{10}{l}{(i) $\alpha = 95\%$}\\
 G1 & 96.60 ( 6.02 )  & 96.77 ( 6.72 )  & 96.50 ( 10.31 )  & 87.97 ( 9.34 )  & 88.07 ( 10.58 )  & 88.4 ( 22.55 )  & 90.33 ( 3.41 )  & 93.63 ( 3.71 )  & 95.07 ( 7.84 )  \\ 
  G2 & 96.17 ( 3.64 )  & 96.67 ( 5.67 )  & 95.43 ( 5.73 )  & 88.8 ( 5.86 )  & 88.80 ( 9.64 )  & 88.93 ( 13.23 )  & 93.73 ( 2.43 )  & 97.57 ( 3.38 )  & 95.00 ( 3.04 )  \\ 
  G3 & 96.07 ( 3.38 )  & 96.30 ( 5.53 )  & 95.57 ( 5.13 )  & 89.13 ( 5.48 )  & 89.00 ( 9.53 )  & 89.8 ( 12.69 )  & 93.97 ( 2.31 )  & 98.00 ( 3.34 )  & 95.67 ( 2.77 )  \\ 
  G4 & 95.93 ( 3.09 )  & 96.63 ( 5.41 )  & 95.40 ( 4.59 )  & 88.53 ( 5.03 )  & 88.30 ( 9.42 )  & 88.67 ( 11.26 )  & 95.33 ( 2.16 )  & 98.03 ( 3.28 )  & 96.00 ( 2.48 )  \\ 
  G5 & 96.17 ( 2.22 )  & 96.4 ( 4.95 )  & 95.47 ( 2.99 )  & 88.83 ( 3.63 )  & 88.63 ( 8.76 )  & 89.07 ( 8.17 )  & 95.77 ( 1.71 )  & 97.77 ( 3.16 )  & 94.90 ( 1.75 )  \\
   \bottomrule
\end{tabular}
\label{tab1:se15}
\end{adjustbox}
\end{table}

\section{
Real Data Analysis
}\label{sec:6}
In this section, we use data from the \href{https://www.census.gov/srd/csrmreports/byyear.html.}{1989 Small Area Income and Poverty Estimates (SAIPE) program} to demonstrate the proposed parametric bootstrap methods for constructing prediction intervals. The SAIPE program produces poverty statistics for various age groups at different levels of geographic aggregation for the 50 states and the District of Columbia. A general area-level model is applied to direct poverty estimates from the U.S. Current Population Survey (CPS), borrowing information from administrative records such as selected tax data from the IRS and participation data from the Supplemental Nutrition Assistance Program (SNAP, formerly known as the Food Stamp Program).

For the 5–17 age group, \cite{bell2006using} detected that the 1989 CPS direct estimate for the state of Connecticut could be regarded as a potential outlier. To address potential outliers, they considered area-level models with $t$-distributed random effects. One such model is given by: 
\begin{eqnarray*} y_i = \bm{x}_i^\prime\bm{\beta} + u_i + e_i, \quad i = 1, \ldots, 51, 
\end{eqnarray*}
where $u_i \overset{iid}{\sim} t_{\phi}(0, A)$ with 3 degrees of freedom (i.e., $\phi = 3$) and $e_i \overset{ind}{\sim} N(0, D_i)$. The SAIPE data provide the direct CPS estimated poverty rates for related children aged 5–17 ($y_i$), associated sampling variances ($D_i$), and four auxiliary variables ($\bm{x}_i$), including pseudo-poverty rates tabulated from IRS tax data, tax non-filer rates from IRS data, Food Stamp participation proportions, and previous census estimated poverty ratios.

By replacing the unknown parameters $\bm{\beta}$ and $A$ with the weighted least squares (WLS) and Fay-Herriot (FH) estimates, respectively, we can obtain the EBLUP, $\hat{\theta}_i$, of the poverty ratios for each state. Following \cite{bell2007use}, we set $\alpha = 90\%$. We then implement both the proposed single ($B_1=400$; \textbf{t,SB}) and double ($B_2 = 100$; \textbf{t,DB}) parametric bootstrap methods to obtain the 90 percent prediction intervals and compare them with those obtained the proposed single parametric bootstrap method under the case $k = \infty$ (i.e., $u_i$ follows a normal distribution; \textbf{normal}) and with direct confidence intervals (\textbf{DIRECT}). The computation times are 32.26 seconds and 136.63 seconds for the proposed single and double bootstrap methods, respectively.

\begin{figure}[!t]
    \centering
    \includegraphics[width=1\linewidth]{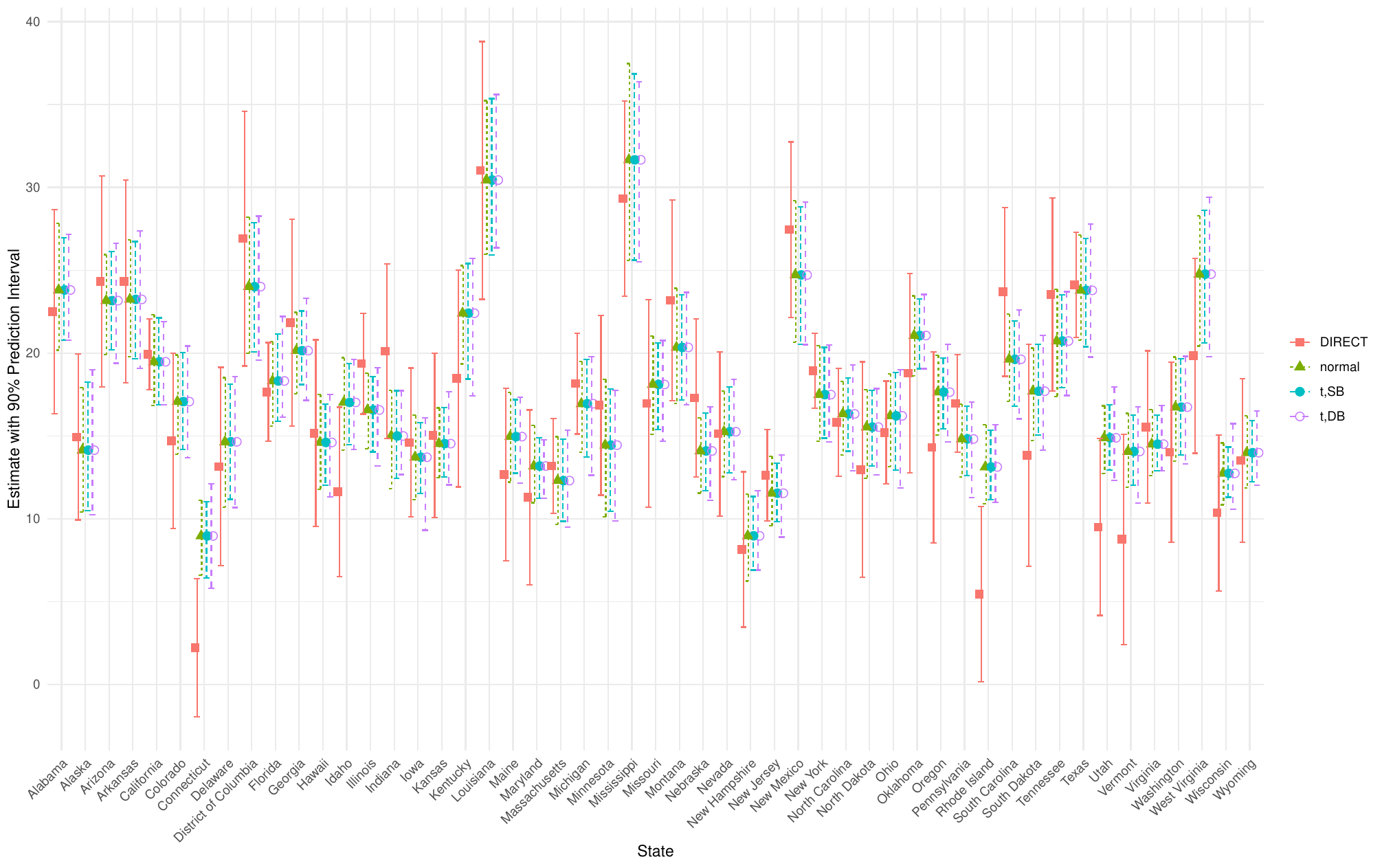}
    \caption{90\% Prediction intervals and estimates for State 5-17 1989 poverty ratios.}
    \label{fig:saipe}
\end{figure}

Figure \ref{fig:saipe} presents the poverty ratio estimates and associated 90\% prediction intervals for each state. For most states, the direct intervals are too wide to be informative. The lengths of the single bootstrap intervals based on normally distributed and $t$-distributed random effects are similar. The double bootstrap intervals are longer than the single bootstrap intervals and tend to contain them. This corresponds with the theory above that double bootstrap intervals could potentially offer better coverage than single bootstrap intervals.

\section{Discussion}\label{sec:7}
In this study, we put forward parametric bootstrap approaches to construct prediction intervals in the contexts of small area estimation under a general (non normal) area level model. Our simulation results show that the proposed single bootstrap method with the Fay-Herriot method of variance estimator performs well for all cases. Moreover, it is more efficient in terms of average lengths than the existing parametric bootstrap method proposed by \cite{hall2006parametric}. 

When the number of area is small, there is a high likelihood to obtain negative estimates of the variance $A$ when applying Prasad Rao method. Throughout our simulation studies, we notice that the likelihood of obtaining negative estimates of variance might affect the performance of the parametric bootstrap methods in terms of both coverage probabilities and average lengths. One might consider a better variance estimator, such as the Fay-Herriot estimator we used, or a sensible truncation of negative estimates. In this study, we arbitrarily truncate the negative estimates of the variance $A$ at 0.01 as similar to \cite{datta2005measuring}. To this end, in the future we will consider an extension of adjusted maximum likelihood estimators of $A$ such as the ones considered by \cite{li2010adjusted} or \cite{hirose2018estimating}  to the model proposed in this paper.

There is an another issue that even though the double bootstrap calibration can bring the coverage accuracy to $o(m^{-1})$ regardless the existence of pivot, our simulations suggest that it is not always beneficial to attempt boosting the theoretical coverage probability via double bootstrap, disregarding other properties of the interval. Specifically, variability of calibrated intervals are greater than uncalibrated ones, minimum length property is almost never preserved, and the results are quite dependent on the parameters and fixed constants of the problem, such as estimation of the variance components $A$ in this study.
For instance, Table \ref{tab1:se50} shows that the proposed single bootstrap intervals already produced good performance and thus the calibration using double bootstrap has little or no effect, where $m = 50$. When $m$ is relatively small, i.e., $m = 15$, double bootstrap improves the coverage probability marginally but it produces much larger interval length than the corresponding single bootstrap method. The substantial increase in interval length may be attributed to numerical instability in the variance estimation when using second-stage bootstrap sample data.

\section*{Disclosure Statement}
The authors report there are no competing interests to declare.

\newpage
\begin{center}
{\large\bf SUPPLEMENTARY MATERIAL}
\end{center}

\textbf{Appendices to Empirical Best Linear Prediction Interval:} These supplementary appendices provide additional technical details and simulation tables for the article. The content includes proofs for Theorems 1 and 2, proofs for Propositions 1, 2, and 3 
in Appendices A-F,
as well as simulation tables that were not included in the main text 
in Appendix G.




\bibliographystyle{chicago}
\bibliography{ref.bib}

@article{bell2007use,
  title={Use of ACS data to produce SAIPE model-based estimates of poverty for counties},
  author={Bell, William and Basel, Wesley and Cruse, Craig and Dalzell, Lucinda and Maples, Jerry and O’Hara, Brett and Powers, David},
  journal={Census Report},
  year={2007},
  publisher={Citeseer}
}

@article{jiang2022goodness,
  title={Goodness-of-fit test with a robustness feature},
  author={Jiang, Jiming and Torabi, Mahmoud},
  journal={TEST},
  volume={31},
  number={1},
  pages={76--100},
  year={2022},
  publisher={Springer}
}

@book{jiang2007linear,
  title={Linear and generalized linear mixed models and their applications},
  author={Jiang, Jiming and Nguyen, Thuan},
  volume={1},
  year={2007},
  publisher={Springer}
}

@article{datta1995robust,
  title={Robust hierarchical Bayes estimation of small area characteristics in the presence of covariates and outliers},
  author={Datta, Gauri S and Lahiri, Partha},
  journal={Journal of Multivariate Analysis},
  volume={54},
  number={2},
  pages={310--328},
  year={1995},
  publisher={Elsevier}
}

@article{hawala2018variance,
  title={Variance modeling for domains},
  author={Hawala, Sam and Lahiri, Partha},
  journal={Statistics and Applications},
  volume={16},
  number={1},
  pages={399--409},
  year={2018}
}

@article{ha2014methods,
  title={Methods and results for small area estimation using smoking data from the 2008 National Health Interview Survey},
  author={Ha, Neung Soo and Lahiri, Partha and Parsons, Van},
  journal={Statistics in Medicine},
  volume={33},
  number={22},
  pages={3932--3945},
  year={2014},
  publisher={Wiley Online Library}
}

@article{mccullough1998implementing,
  title={Implementing the double bootstrap},
  author={McCullough, BD and Vinod, HD},
  journal={Computational Economics},
  volume={12},
  pages={79--95},
  year={1998},
  publisher={Springer}
}

@article{fabrizi2010robust,
  title={Robust linear mixed models for small area estimation},
  author={Fabrizi, Enrico and Trivisano, Carlo},
  journal={Journal of Statistical Planning and Inference},
  volume={140},
  number={2},
  pages={433--443},
  year={2010},
  publisher={Elsevier}
}

@book{hall2013bootstrap,
  title={The bootstrap and Edgeworth expansion},
  author={Hall, Peter},
  year={2013},
  pages = {17},
  publisher={Springer Science \& Business Media}
}

@book{rao2015small,
  title={Small area estimation},
  author={Rao, John NK and Molina, Isabel},
  year={2015},
  publisher={John Wiley \& Sons}
}

@article{cox1975prediction,
  title={Prediction intervals and empirical Bayes confidence intervals},
  author={Cox, DR},
  journal={Journal of Applied Probability},
  volume={12},
  number={S1},
  pages={47--55},
  year={1975},
  publisher={Cambridge University Press}
}

@article{yoshimori2014second,
author = {Masayo Yoshimori and Partha Lahiri},
title = {{A second-order efficient empirical Bayes confidence interval}},
volume = {42},
journal = {The Annals of Statistics},
number = {4},
publisher = {Institute of Mathematical Statistics},
pages = {1233 -- 1261},
year = {2014},
doi = {10.1214/14-AOS1219},
URL = {https://doi.org/10.1214/14-AOS1219}
}

@inproceedings{bell2006using,
  title={Using the t-distribution to deal with outliers in small area estimation},
  author={Bell, William R and Huang, Elizabeth T},
  booktitle={Proceedings of statistics Canada symposium},
  year={2006}
}

@article{xie2007estimation,
  title={Estimation of the proportion of overweight individuals in small areas—a robust extension of the Fay--Herriot model},
  author={Xie, Dawei and Raghunathan, Trivellore E and Lepkowski, James M},
  journal={Statistics in Medicine},
  volume={26},
  number={13},
  pages={2699--2715},
  year={2007},
  publisher={Wiley Online Library}
}

@article{chatterjee2008parametric,
author = {Snigdhansu Chatterjee and Partha Lahiri and Huilin Li},
title = {{Parametric bootstrap approximation to the distribution of EBLUP and related prediction intervals in linear mixed models}},
volume = {36},
journal = {The Annals of Statistics},
number = {3},
publisher = {Institute of Mathematical Statistics},
pages = {1221 -- 1245},
year = {2008},
doi = {10.1214/07-AOS512},
URL = {https://doi.org/10.1214/07-AOS512}
}

@article{hall2006parametric,
  title={On parametric bootstrap methods for small area prediction},
  author={Hall, Peter and Maiti, Tapabrata},
  journal={Journal of the Royal Statistical Society Series B: Statistical Methodology},
  volume={68},
  number={2},
  pages={221--238},
  year={2006},
  publisher={Oxford University Press}
}

@article{shi1992accurate,
  title={Accurate and efficient double-bootstrap confidence limit method},
  author={Shi, Sheng G},
  journal={Computational statistics \& data analysis},
  volume={13},
  number={1},
  pages={21--32},
  year={1992},
  publisher={Elsevier}
}

@article{li2010adjusted,
  title={An adjusted maximum likelihood method for solving small area estimation problems},
  author={Li, Huilin and Lahiri, Partha},
  journal={Journal of multivariate analysis},
  volume={101},
  number={4},
  pages={882--892},
  year={2010},
  publisher={Elsevier}
}

@article{hirose2018estimating,
  title={Estimating variance of random effects to solve multiple problems simultaneously},
  author={Hirose, Masayo Yoshimori and Lahiri, Partha},
  journal={The Annals of Statistics},
  volume={46},
  number={4},
  pages={1721--1741},
  year={2018},
  publisher={JSTOR}
}

@article{ghosh1987robust,
  title={Robust empirical Bayes estimation of means from stratified samples},
  author={Ghosh, Malay and Lahiri, Parthasarathi},
  journal={Journal of the American Statistical Association},
  volume={82},
  number={400},
  pages={1153--1162},
  year={1987},
  publisher={Taylor \& Francis}
}

@article{fay1979estimates,
  title={Estimates of income for small places: an application of James-Stein procedures to census data},
  author={Fay, Robert E and Herriot, Roger A},
  journal={Journal of the American Statistical Association},
  volume={74},
  number={366a},
  pages={269--277},
  year={1979},
  publisher={Taylor \& Francis}
}

@article{prasad1990estimation,
  title={The estimation of the mean squared error of small-area estimators},
  author={Prasad, NG Narasimha and Rao, Jon NK},
  journal={Journal of the American statistical association},
  volume={85},
  number={409},
  pages={163--171},
  year={1990},
  publisher={Taylor \& Francis}
}

@article{datta2005measuring,
  title={On measuring the variability of small area estimators under a basic area level model},
  author={Datta, Gauri Sankar and Rao, JNK and Smith, David Daniel},
  journal={Biometrika},
  volume={92},
  number={1},
  pages={183--196},
  year={2005},
  publisher={Oxford University Press}
}

@book{shao2008mathematical,
  title={Mathematical statistics},
  author={Shao, Jun},
  year={2008},
  pages={471–-477},
  publisher={Springer Science \& Business Media}
}

@article{lahiri1995robust,
  title={Robust estimation of mean squared error of small area estimators},
  author={Lahiri, P and Rao, JNK},
  journal={Journal of the American Statistical Association},
  volume={90},
  number={430},
  pages={758--766},
  year={1995},
  publisher={Taylor \& Francis}
}

@inproceedings{otto1995sampling,
  title={Sampling error modelling of poverty and income statistics for states},
  author={Otto, Mark C and Bell, William R},
  booktitle={American Statistical Association, Proceedings of the Section on Government Statistics},
  pages={160--165},
  year={1995}
}
\end{document}